\documentclass[10pt]{article}
\usepackage{color}
\usepackage{amssymb,amsmath,graphicx,acronym,amsfonts}
\usepackage[square, comma, sort&compress, numbers]{natbib}

\newtheorem{corollary}{Corollary}[section]
\newtheorem{theorem}{Theorem}[section]
\newtheorem{lemma}{Lemma}[section]
\newtheorem{definition}{Definition}[section]
\newtheorem{proposition}{Proposition}[section]
\newtheorem{example}{Example}[section]
\newtheorem{assum}{Assumption}[section]
\newtheorem{algo}{Algorithm}[section]
\newtheorem{Remark}{Remark}[section]

\def\bc{\begin{corl}}
\def\bc{\end{corl}}
\def\ba{\begin{algo}}
\def\ea{\end{algo}}
\def\br{\begin{Remark}}
\def\er{\end{Remark}}
\def\bs{\begin{assum}}
\def\es{\end{assum}}
\def\bt{\begin{theorem}}
\def\et{\end{theorem}\vskip 3pt}
\def\bl{\begin{lemma}}
\def\el{\end{lemma}}
\def\ep{\end{proposition}}
\def\bp{\begin{proposition}}
\def\qed{\hfill{$\Box$}\vskip 5pt}
\def\be{\begin{example}}
\def\ee{\end{example}}
\def\bd{\begin{definition}}
\def\ed{\end{definition}}
\def\bc{\begin{corollary}}
\def\ec{\end{corollary}}
\def\proof{\noindent\it Proof. \hspace{1mm}\rm}
\input epsf
\begin{document}
\title{\bf \Large Copositivity Detection of Tensors: Theory and Algorithm}
\author{Haibin Chen\thanks{School of Management Science, Qufu Normal University, Rizhao, Shandong, P.R. China,
Department of Applied Mathematics, The Hong Kong Polytechnic University, Hung Hom,
Kowloon, Hong Kong. Email: chenhaibin508@163.com. } \quad
Zheng-Hai Huang
\thanks{Department of Mathematics, School of Science, Tianjin University, Tianjin 300072, P.R. China. %This author is also with the Center for Applied Mathematics of Tianjin University.
Email: huangzhenghai@tju.edu.cn. This author's work was supported by the National Natural Science Foundation of China (Grant No. 11431002).} \quad
Liqun Qi
\thanks{Department of Applied Mathematics, The Hong Kong Polytechnic University, Hung Hom,
Kowloon, Hong Kong. Email: maqilq@polyu.edu.hk. This author's work
was supported by the Hong Kong Research Grant Council (Grant No.
PolyU 502111, 501212, 501913 and 15302114).}}
\date{}
\maketitle
\vspace{-0.6cm}
\begin{abstract}
A symmetric tensor is called copositive if it generates a multivariate form taking nonnegative values over the nonnegative orthant. Copositive tensors have found important applications in polynomial optimization and tensor complementarity problems. In this paper, we consider copositivity detection of tensors both from theoretical and computational points of view. After giving several necessary conditions for copositive tensors,
we propose several new criteria for copositive tensors based on the representation of the multivariate form in barycentric coordinates with respect to the standard simplex and simplicial partitions. It is verified that, as the partition gets finer and finer, the concerned conditions eventually capture all strictly copositive tensors. Based on the obtained theoretical results with the help of simplicial partitions, we propose a numerical method to judge whether a tensor is copositive or not. The preliminary numerical results confirm our theoretical findings.
\medskip



\end{abstract}

\newpage
\section{Introduction}
A symmetric tensor is called copositive if it generates a multivariate form taking nonnegative values over the nonnegative orthant\cite{qlq2013}.
Copositive tensors constitute a large class of tensors
that contain nonnegative tensors and several kinds of structured tensors in the even order symmetric case such as $M$-tensors, diagonally dominant tensors and so on\cite{chen14,Ding15,Kannan15,LCQL,LCL,LWZ,Qi14,qi14,Zhang12}. Recently, it has been found that copositive tensors have important applications in polynomial optimization \cite{Pena14,song2016} and the tensor complementarity problem \cite{CQW, SQ15, SQ}.

Pena et al. \cite{Pena14} provided a general characterization for a class of polynomial optimization problems that can be formulated as a conic program over the cone of completely positive tensors or copositive tensors. As a consequence of this characterization, it follows that recent related results for quadratic problems can be further strengthened and generalized to higher order polynomial optimization problems. On the other hand, Che, Qi and Wei \cite{CQW} showed that the tensor complementarity problem with a strictly copositive tensor has a nonempty and compact solution set; Song and Qi \cite{SQ15} proved that a real symmetric tensor is a (strictly) semi-positive  if and only if it is (strictly) copositive, and Song and Qi \cite{SQ15,SQ} obtained several results for the tensor complementarity problem with a (strictly) semi-positive tensor; and Huang and Qi \cite{HQ2016} formulated an $n$-person noncooperative game as a tensor complementarity problem with the involved tensor is nonnegative. Thus, copositive tensors play an important role in the tensor complementarity problem.
Now, there is a challenging problem that is how to check the copositivity of a given symmetric tensors efficiently?

Though many structured tensors are copositive, from the practical point of view, it is of great important to check directly whether a tensor is copositive or not. Several sufficient conditions or necessary and sufficient conditions for copositive tensors have been presented in \cite{qlq2013,Song15}. However, it is hard to verify numerically whether a tensor is copositive or not from these conditions. In \cite{song2016}, Song and Qi gave the concepts of Pareto H-eigenvalue and Pareto Z-eigenvalue for symmetric tensors. It is proved that a symmetric tensor $\mathcal{A}$ is strictly copositive if and only if every
Pareto H-eigenvalue (Z-eigenvalue) of $\mathcal{A}$ is positive, and $\mathcal{A}$ is copositive if and only if every
Pareto H-eigenvalue (Z-eigenvalue) of $\mathcal{A}$ is nonnegative. Unfortunately, it is NP-hard the compute the minimum Pareto H-eigenvalue or Pareto Z-eigenvalue of a given symmetric tensors.
In fact, for a given tensor, the problem to judge whether it is copositive or not is NP-complete, even for the matrix case \cite{Dickinson14,Murty87}. To the best of our knowledge, there is not any numerical detection method for copositive tensors with order greater than three. In this paper, we further give some theoretical studies on various conditions for (strictly) copositive tensors; and based on some of our theoretical findings, we propose a numerical method to judge whether a tensor is copositive or not. The algorithm we investigated can be viewed as an extension of some branch-and-bound type algorithm for testing copositivity of symmetric matrices \cite{Bundfuss08,Sponsel12,Xu11}.

The rest of this paper is organized as follows. In Section 2, we recall some notions and basic facts about tensors and the corresponding homogeneous polynomials. Three necessary conditions for copositive tensors are given in Section 3.  In Section 4, we give several criteria for (strictly) copositive tensors based on the simplicial subdivision, and an equivalent condition for a symmetric tensor that is not copositive. In Section 5, we propose a numerical detection algorithm for copositive tensors based on the results obtained in Section 4; and show that the algorithm can always capture strictly copositive tensors in finitely many iterations. The preliminary numerical results are reported in Section 6, and final remarks and some future work are given in Section 7.

\setcounter{equation}{0}
\section{Preliminaries}

Throughout this paper, we denote the set consisting of all positive integers by $\mathbb{N}$, and always assume that $m, n\in \mathbb{N}$. Let $\mathbb{R}^n$ be the $n$ dimensional real Euclidean space and the set of all nonnegative vectors in $\mathbb{R}^n$ be denoted by $\mathbb{R}^n_+$. Let $\mathbb{R}^n_{++}$ denote the set of vectors with positive entries. Vectors are denoted by bold lowercase letters i.e. ${\bf x},~ {\bf y},\cdots$, matrices are denoted by capital letters i.e. $A, B, \cdots$, and tensors are written as calligraphic capitals such as $\mathcal{A}, \mathcal{T}, \cdots$. We denote $[n]=\{1,2,\cdots,n\}$. The $i$-th unit coordinate vector in $\mathbb{R}^n$ is denoted by ${\bf e_i}$ for any $i\in [n]$.

A real $m$-th order $n$-dimensional tensor $\mathcal{A}=(a_{i_1i_2\cdots i_m})$ is a multi-array of real entries $a_{i_1i_2\cdots i_m}$, where $i_j \in [n]$ for $j\in [m]$. In this paper, we always assume that $m\geq 3$ and $n\geq 2$. A tensor is said to be nonnegative if all its entries are nonnegative. If the entries $a_{i_1i_2\cdots i_m}$ are invariant under any permutation of their indices, then tensor $\mathcal{A}$ is called a symmetric tensor. In this paper, we always consider real symmetric tensors. The identity tensor $\mathcal{I}$ with order $m$ and dimension $n$ is given by $\mathcal{I}_{i_1\cdots i_m}=1$ if $i_1=\cdots=i_m$ and $\mathcal{I}_{i_1\cdots i_m}=0$ otherwise. For any $J\subseteq [n]$, $|J|$ denotes the number of elements of $J$, and ${\cal A}_J$ denotes a principle subtensor of $\mathcal{A}$.

We denote
$$
S_{m,n}:=\{\mathcal{A}: \mathcal{A} \mbox{ is an } m\mbox{-th~order }  n\mbox{-dimensional} \mbox{  symmetric tensor}\}.
$$
Clearly, $S_{m,n}$ is a vector space under the addition and multiplication defined as below: for any $t \in \mathbb{R}$,
$\mathcal{A}=(a_{i_1 \cdots i_m})_{1 \le i_1,\cdots,i_m
\le n}$ and $\mathcal{B}=(b_{i_1 \cdots i_m})_{1 \le
i_1,\cdots,i_m \le n},$
\[
\mathcal{A}+\mathcal{B}=(a_{i_1 \cdots i_m}+b_{i_1 \cdots i_m})_{1
\le i_1,\cdots,i_m \le n}\quad \mbox{\rm and }\quad t
\mathcal{A}=(ta_{i_1 \cdots i_m})_{1 \le i_1,\cdots,i_m
\le n}.
\]

For any $\mathcal{A}, \mathcal{B} \in S_{m,n}$, we define the inner product by
$\langle \mathcal{A},\mathcal{B}\rangle:=\sum_{i_1,\cdots,i_m=1}^{n}a_{i_1 \cdots i_m}b_{i_1 \cdots i_m}$, and the corresponding norm by
$$
\|\mathcal{A}\|=\left(\langle \mathcal{A},\mathcal{A}\rangle\right)^{1/2}=\left(\sum_{i_1,\cdots,i_m=1}^{n}(a_{i_1
\cdots i_m})^2\right)^{1/2}.
$$

For any ${\bf x}\in \mathbb{R}^n$, we use $x_i$ to denote its $i$th component; and use $\|{\bf x}\|_m$ to denote the $m$-norm of ${\bf x}$.

For $m$ vectors ${\bf x},{\bf y}, \cdots, {\bf z}\in
\mathbb{R}^n$, we use ${\bf x}\circ {\bf y}\circ \cdots \circ {\bf z}$ to denote the $m$-th order $n$-dimensional symmetric rank one tensor with
\[
({\bf x}\circ {\bf y}\circ \cdots \circ {\bf z})_{i_1 i_2\cdots i_m}=x_{i_1}y_{i_2}\cdots z_{i_m}, \ \forall \, i_1,\cdots,i_m \in [n].
\]
And the inner product of a symmetric tensor and the rank one tensor is given by
$$
\langle
\mathcal{A}, {\bf x}\circ {\bf y}\circ \cdots \circ {\bf z}\rangle:=\sum_{i_1,\cdots,i_m=1}^{n}a_{i_1 \cdots
i_m}x_{i_1}y_{i_2}\cdots z_{i_m}.
$$
For $m\in \mathbb{N}$ and $k\in [m]$, we denote
$$
\mathcal{A}{\bf x}^k{\bf y}^{m-k}=\langle\mathcal{A}, \underbrace{{\bf x}\circ \cdots {\bf x}}_k\circ\underbrace{{\bf y}\circ \cdots \circ {\bf y}}_{m-k}\rangle\quad \mbox{\rm and}\quad \mathcal{A}{\bf x}^m=\langle\mathcal{A}, \underbrace{{\bf x}\circ \cdots {\bf x}}_m\rangle,
$$
then
\begin{equation}\label{e21}
\mathcal{A}{\bf x}^k{\bf y}^{m-k}=\sum_{i_1,\cdots,i_m=1}^{n}a_{i_1 \cdots i_m}x_{i_1}\cdots x_{i_k}y_{i_{k+1}}\cdots y_{i_m}
\quad \mbox{\rm and}\quad \mathcal{A}{\bf x}^m=\sum_{i_1,\cdots,i_m=1}^{n}a_{i_1 \cdots i_m}x_{i_1}\cdots x_{i_m}.
\end{equation}
For any $\mathcal{A}=(a_{i_1i_2\cdots i_m})\in S_{m,n}$ and ${\bf x}\in \mathbb{R}^n$, we have $\mathcal{A}{\bf x}^{m-1}\in \mathbb{R}^n$ with
$$
(\mathcal{A}{\bf x}^{m-1})_i=\sum_{i_2,i_3,\cdots,i_m\in [n]}a_{ii_2\cdots i_m}x_{i_2}\cdots x_{i_m},~~\forall \, i\in [n].
$$

It is known that an $m$-th order $n$-dimensional symmetric tensor defines uniquely an $m$-th degree homogeneous
polynomial $f_{\mathcal{A}}({\bf x})$ on $\mathbb{R}^n$: for all ${\bf x}=(x_1,\cdots,x_n)^T
\in \mathbb{R}^n$,
$$
f_{\mathcal{A}}({\bf x})= \mathcal{A}{\bf x}^m=\sum_{i_1,i_2,\cdots, i_m\in [n]}a_{i_1i_2\cdots i_m}x_{i_1}x_{i_2}\cdots x_{i_m};
$$
and conversely, any $m$-th degree homogeneous polynomial function $f({\bf x})$ on $\mathbb{R}^n$ also corresponds uniquely a symmetric tensor. Furthermore, an even order symmetric tensor $\mathcal{A}$ is called positive semi-definite (positive definite) if $f_{\mathcal{A}}({\bf x}) \geq 0$ ($f_{\mathcal{A}}({\bf x})> 0$) for all ${\bf x}\in \mathbb{R}^n$ (${\bf x}\in \mathbb{R}^n \backslash \{\bf 0\}$).

\setcounter{equation}{0}
\section{Necessary conditions for copositivity of tensors}
In this section, we first introduce the definition of copositive tensors, and then three necessary conditions are established, which are given based on binomial expansion, principle subtensor and convex combination of symmetric tensors, respectively.

\bd\label{def31} \cite{qlq2013} Let $\mathcal{A}\in S_{m,n}$ be given. If
$\mathcal{A}{\bf x}^m\geq 0\;(\mathcal{A}{\bf x}^m>0)$ for any ${\bf x}\in \mathbb{R}^n_+({\bf x}\in \mathbb{R}^n_+\backslash \{{\bf0}\})$, then
$\mathcal{A}$ is called a copositive (strictly copositive) tensor.
\ed

To move on, we first prove the following result, which looks like the binomial expansion of two real variables.
\bl\label{lema31} Let $\mathcal{A}=(a_{i_1i_2\cdots i_m})\in S_{m,n}$ be given. For any ${\bf x}, {\bf y}\in \mathbb{R}^n$, it holds that
$$
\mathcal{A}({\bf x}+{\bf y})^m=\mathcal{A}{\bf x}^m+\binom m 1 \mathcal{A}{\bf x}^{m-1}{\bf y}+\binom m 2\mathcal{A}{\bf x}^{m-2}{\bf y}+\cdots+\binom m m\mathcal{A}{\bf y}^m.
$$
\el
\proof  For any ${\bf x}, {\bf y}\in \mathbb{R}^n$, we know that
\begin{equation*}%\label{e23}
\begin{aligned}
\mathcal{A}({\bf x}+{\bf y})^m=&\langle\mathcal{A},({\bf x}+{\bf y})\circ ({\bf x}+{\bf y})\circ\cdots \circ ({\bf x}+{\bf y})\rangle \\
=&\sum_{i_1,i_2,\cdots, i_m\in [n]}a_{i_1i_2\cdots i_m}(x_{i_1}+y_{i_1})(x_{i_2}+y_{i_2}) \cdots (x_{i_m}+y_{i_m}).
\end{aligned}
\end{equation*}
By using the symmetry property of $\mathcal{A}$ and (\ref{e21}), we further obtain that
$$
\begin{aligned}
\mathcal{A}({\bf x}+{\bf y})^m=&\sum_{i_1,i_2,\cdots, i_m\in [n]}a_{i_1i_2\cdots i_m}x_{i_1}x_{i_2}\cdots x_{i_m} \\
&+\binom m 1\sum_{i_1,i_2,\cdots, i_m\in [n]}a_{i_1i_2\cdots i_m}x_{i_1}\cdots x_{i_{m-1}}y_{i_m} \\
&+\binom m 2\sum_{i_1,i_2,\cdots, i_m\in [n]}a_{i_1i_2\cdots i_m}x_{i_1}\cdots x_{i_{m-2}}y_{i_{m-1}}y_{i_m}\\
&+\cdots \cdots \\
&+\binom m m\sum_{i_1,i_2,\cdots, i_m\in [n]}a_{i_1i_2\cdots i_m}y_{i_1}y_{i_2}\cdots y_{i_m} \\
=&\mathcal{A}{\bf x}^m+\binom m 1 \mathcal{A}{\bf x}^{m-1}{\bf y}+\binom m 2\mathcal{A}{\bf x}^{m-2}{\bf y}+\cdots+\binom m m\mathcal{A}{\bf y}^m,
\end{aligned}
$$
which completes the proof.
\qed

Now, we use Lemma \ref{lema31} to establish a necessary condition for the copositive tensor.

\bt\label{them31} Let $\mathcal{A}\in S_{m,n}$ be copositive. If there is ${\bf x}\in \mathbb{R}^n_+$ such that $\mathcal{A}{\bf x}^m=0$, then $\mathcal{A}{\bf x}^{m-1}\geq {\bf 0}$.
\et
\proof Suppose ${\bf u}\in \mathbb{R}^n_+$ and $\varepsilon>0, \varepsilon\in \mathbb{R}$. Since $\mathcal{A}$ is copositive, it follows that
$$
\mathcal{A}{\bf z}^m\geq 0,~~\mbox{\rm where}~~{\bf z}={\bf x}+\varepsilon {\bf u}.
$$
So, by Lemma \ref{lema31}, we have
\begin{equation}\label{e31}
\begin{aligned}
\mathcal{A}{\bf z}^m=& \mathcal{A}({\bf x}+\varepsilon{\bf u})^m \\
=&\mathcal{A}{\bf x}^m+\binom m 1 \varepsilon \mathcal{A}{\bf x}^{m-1}{\bf u}+\binom m 2 \varepsilon^2 \mathcal{A}{\bf x}^{m-2}{\bf u}^2+\cdots+\binom m m\varepsilon^m\mathcal{A}{\bf u}^m.
\end{aligned}
\end{equation}
Now, assume that there is an index $i\in [n]$ such that $(\mathcal{A}{\bf x}^{m-1})_i<0$. Take ${\bf u}={\bf e_i}$ in (\ref{e31}), we obtain that
$$
\begin{aligned}
\mathcal{A}{\bf z}^m=& \mathcal{A}({\bf x}+\varepsilon{\bf e_i})^m \\
=&\mathcal{A}{\bf x}^m+\binom m 1 \varepsilon \mathcal{A}{\bf x}^{m-1}{\bf e_i}+\cdots+\binom m m\varepsilon^m\mathcal{A}{\bf e_i}^m \\
=& \binom m 1 \varepsilon (\mathcal{A}{\bf x}^{m-1})_i+\binom m 2 \varepsilon^2 \mathcal{A}{\bf x}^{m-2}{\bf e_i}^2+\cdots+\binom m m\varepsilon^m\mathcal{A}{\bf e_i}^m \\
<&0
\end{aligned}
$$
holds for any sufficiently small $\varepsilon>0$, which contradicts the fact that $\mathcal{A}$ is copositive; and hence, the desired result follows. \qed

Next, we establish a necessary condition for the copositive tensor by using the concept of principle subtensor.

\bt\label{them32} Let $\mathcal{A}\in S_{m,n}$ be copositive. Then, for any $J\subseteq [n]$, the system $\mathcal{A}_J{\bf x}^{m-1}\geq {\bf 0}$ admits a nonzero solution ${\bf x}\in \mathbb{R}_+^{|J|}$.
\et
\proof We prove it by contradiction. If there is $J\subseteq [n]$ such that
\begin{equation}\label{e32}
\mathcal{A}_J{\bf x}^{m-1}< {\bf 0}~~\mbox{\rm for~all~nonzero}~~{\bf x}\in \mathbb{R}_+^{|J|},
\end{equation}
then, we may define a vector ${\bf y}\in \mathbb{R}^n$ with $y_i=\frac{x_i}{\|{\bf x}\|_m}$ when $i\in J$; and $y_i=0$ otherwise. By a direct computation,
we obtain $\mathcal{A}{\bf y}^m=\frac{\mathcal{A}_J{\bf x}^m}{\|{\bf x}\|_m^m}<0$, which contradicts that $\mathcal{A}$ is copositive.
\qed

To end this section, we establish a necessary condition for the copositivity of a convex combination of two symmetric tensors.

\bt\label{them33} Let $\mathcal{A},\mathcal{B}\in S_{m,n}$ be given.
Suppose that there exists $t\in [0,1]$ such that $(1-t)\mathcal{A}+t\mathcal{B}$ is copositive, then
$\max \{\mathcal{A}{\bf u}^m+\mathcal{A}{\bf v}^m, \mathcal{B}{\bf u}^m+\mathcal{B}{\bf v}^m\}\geq 0$ for all ${\bf u}, {\bf v}\in \mathbb{R}^n_+$.
\et
\proof Since $(1-t)\mathcal{A}+t\mathcal{B}$ is copositive, it follows that
$$
(1-t)\mathcal{A}{\bf u}^m+t\mathcal{B}{\bf u}^m\geq 0,~~(1-t)\mathcal{A}{\bf v}^m+t\mathcal{B}{\bf v}^m\geq 0
$$
for all ${\bf u}, {\bf v}\in \mathbb{R}^n_+$. Adding these two inequalities, one has that
$$
\max \{\mathcal{A}{\bf u}^m+\mathcal{A}{\bf v}^m, \mathcal{B}{\bf u}^m+\mathcal{B}{\bf v}^m\}\geq(1-t)(\mathcal{A}{\bf u}^m+\mathcal{A}{\bf v}^m)+t(\mathcal{B}{\bf u}^m+\mathcal{B}{\bf v}^m)\geq 0
$$
holds for all ${\bf u}, {\bf v}\in \mathbb{R}^n_+$. This completes the proof.
\qed
When $m=2$, Crouzeix, Martin$\acute{e}$z-Legaz, and Seeger proved the conclusions of Theorem \ref{them33} is sufficient and necessary in the matrix case \cite{Crou95}.

\setcounter{equation}{0}
\section{Detection criteria based on simplicial partition}
In this section, several sufficient conditions or necessary conditions of copositive tensors are characterized based on some simplices.
First of all, we show an useful result obtained by Song and Qi \cite{Song15}.
\bl\label{lema41} Let $\mathcal{A}\in S_{m,n}$ be given and $\|\cdot\|$ denote any norm on $\mathbb{R}^n$. Then, we have

{\bf (i)} $\mathcal{A}$ is copositive if and only if $\mathcal{A}{\bf x}^m\geq 0$ for all ${\bf x}\in \mathbb{R}^n_+$ with $\|{\bf x}\|=1$;

{\bf (ii)} $\mathcal{A}$ is strictly copositive if and only if $\mathcal{A}{\bf x}^m> 0$ for all ${\bf x}\in \mathbb{R}^n_+$ with $\|{\bf x}\|=1$.
\el

It is well known that the set $S_0=\{{\bf x}\in \mathbb{R}^n_+~|~\|{\bf x}\|_1=1 \}$ is the so-called standard simplex with vertices ${\bf e_1}, {\bf e_2},\cdots,{\bf e_n}$. So, it follows from Lemma \ref{lema41} that the copositivity of tensor $\mathcal{A}\in S_{m,n}$ can be translated to check
$$
f({\bf x})=\mathcal{A}{\bf x}^m\geq 0,~~\forall~{\bf x}\in S_0.
$$
Thus, our main goal in this section is to search for conditions that can guarantee the homogeneous polynomial $f({\bf x})$ to be nonnegative
on a simplex. A simple way to describe a polynomial with respect to a simplex is to use barycentric coordinates, which gives a convenient
verifiable sufficient condition for a tensor to be copositive on a simplex. This approach has been much used for
convex surface fitting in computer aided geometric design \cite{Juttler00} and the copositivity detection of matrices \cite{Bundfuss08,Sponsel12,Xu11}.

\bl\label{lema42} Let $S_1=conv\{ {\bf u}_1,{\bf u}_2,\cdots,{\bf u}_n\}$ be a simplex. If
\begin{equation}\label{e41}
\langle \mathcal{A},~{\bf u}_{i_1}\circ {\bf u}_{i_2} \cdots \circ {\bf u}_{i_m} \rangle \geq 0\;\left(\langle \mathcal{A},~{\bf u}_{i_1}\circ {\bf u}_{i_2} \cdots \circ {\bf u}_{i_m} \rangle >0\right)
\end{equation}
for all $i_1,i_2,\cdots,i_m \in [n]$, then $\mathcal{A}{\bf x}^m\geq 0\;$($\mathcal{A}{\bf x}^m>0$ respectively) for all ${\bf x}\in S_1$.
\el
\proof For any ${\bf x}\in S_1$, we have that
\begin{equation}\label{e42}
{\bf x}=\lambda_1{\bf u}_1+\cdots+\lambda_n{\bf u}_n,\quad \sum_{i=1}^n\lambda_i=1,\quad \lambda_i\geq 0,\; \forall \, i\in [n].
\end{equation}
So, it holds that
$$
\begin{aligned}
\mathcal{A}{\bf x}^m=&\langle \mathcal{A},~(\lambda_1{\bf u}_1+\cdots+\lambda_n{\bf u}_n)^m\rangle \\
=&\sum_{j_1,j_2,\cdots,j_m\in [n]}\lambda_{j_1}\lambda_{j_2}\cdots \lambda_{j_m}\langle\mathcal{A},~{\bf u}_{j_1}\circ {\bf u}_{j_2}\cdots \circ {\bf u}_{j_m}\rangle.
\end{aligned}
$$
By (\ref{e41}) and (\ref{e42}), we further obtain that
$$
\mathcal{A}{\bf x}^m\geq 0\; ({\bf or}~\mathcal{A}{\bf x}^m> 0 ),
$$
which completes the proof. \qed
If we apply Lemma \ref{lema42} to the standard simplex $S_0=\{{\bf x}\in \mathbb{R}^n_+~|~\|{\bf x}\|_1=1 \}$, it shows that
$$
\langle \mathcal{A},~{\bf e_{i_1}}\circ {\bf e_{i_2}}\circ \cdots \circ {\bf e_{i_m}} \rangle=a_{i_1i_2\cdots i_m} \geq 0, ~\forall~ i_1,i_2,\cdots,i_m \in [n].
$$
This means that all nonnegative tensors are copositive tensors.

Let $S,S_1,S_2,\cdots,S_r$ be finite simplices in $\mathbb{R}^n$. The set $\tilde{S}=\{S_1,S_2,\cdots,S_r\}$ is called a simplicial partition of $S$ if it satisfies that
$$
S=\bigcup_{i=1}^rS_i\quad \mbox{\rm and}\quad \mbox{\rm int} S_i\bigcap \mbox{\rm int} S_j=\emptyset \;\, \mbox{\rm for any}\;\, i,j\in [r]\;\, \mbox{\rm with}\; i\neq j,
$$
where $\mbox{\rm int}S$ denotes the interior of $S$.
\bt\label{theorem41} Let $\mathcal{A}\in S_{m,n}$ be given. Suppose $\widetilde{S}=\{S_1,S_2,\cdots,S_r\}$ is a simplicial partition of simplex $S_0=\{{\bf x}\in \mathbb{R}^n_+~|~\|{\bf x}\|_1=1 \}$; and the vertices of simplex $S_k$ are denoted by ${\bf u}^k_1,{\bf u}^k_2,\cdots,{\bf u}^k_n$ for any $k\in [r]$. Then, the following results hold:

{\bf (i)} if $\langle \mathcal{A},~{\bf u}^k_{i_1}\circ {\bf u}^k_{i_2}\circ \cdots \circ {\bf u}^k_{i_m} \rangle \geq 0$ for all $k\in [r], i_j\in [n],j\in [m]$, then $\mathcal{A}$ is copositive;

{\bf (ii)} if $\langle \mathcal{A},~{\bf u}^k_{i_1}\circ {\bf u}^k_{i_2}\circ \cdots \circ {\bf u}^k_{i_m} \rangle > 0$ for all $k\in [r], i_j\in [n],j\in [m]$, then $\mathcal{A}$ is strictly copositive.
\et
\proof By Lemma \ref{lema41}, it suffices to prove that
$$
\mathcal{A}{\bf x}^m\geq0\;\; \mbox{\rm for~all}~{\bf x}\in S_0.
$$
For any ${\bf x}\in S_0$, since $\widetilde{S}$ is a simplicial partition of $S_0$, it follows that there is an index $k\in [r]$ such that ${\bf x}\in S_k\subseteq \widetilde{S}$. By assumptions and Lemma \ref{lema42}, the desired results follow. \qed

It is easy to see that a simplex $S$ is determined by its vertices, which can be further represented by a matrix $V_S$ whose
columns are vertices of the simplex. It is obvious that $V_S$ is nonsingular and unique up to a permutation of its columns.
So, analogue to Theorem \ref{theorem41}, we have the following results.

\bt\label{theorem42}
Let $\mathcal{A}\in S_{m,n}$ be given. Suppose $\widetilde{S}=\{S_1,S_2,\cdots,S_r\}$ is a simplicial partition of simplex $S_0=\{{\bf x}\in \mathbb{R}^n_+~|~\|{\bf x}\|_1=1 \}$; and the vertices of simplex $S_k$ are denoted by ${\bf u}^k_1,{\bf u}^k_2,\cdots,{\bf u}^k_n$ for any $k\in [r]$. Let $V_{S_k}=({\bf u}^k_1\;{\bf u}^k_2\;\cdots\;{\bf u}^k_n)$ be the matrix corresponding to simplex $S_k$ for any $k\in [r]$. Then, the following results hold:

{\bf (i)} if $V_{S_k}^T\mathcal{A}V_{S_k}$ is copositive for all $k\in [r]$, then $\mathcal{A}$ is copositive;

{\bf (ii)} if $V_{S_k}^T\mathcal{A}V_{S_k}$ is strictly copositive for all $k\in [r]$, then $\mathcal{A}$ is strictly copositive.
\et
\proof It is sufficient to prove {\bf (i)}, since {\bf (ii)} can be verified similarly. For any $k\in [r]$ and $i_1,i_2,\cdots,i_m\in [n]$, we have that
\begin{equation}\label{e43}
\begin{aligned} (V_{S_k}^T\mathcal{A}V_{S_k})_{i_1i_2\cdots i_m}=&\sum_{j_1,j_2,\cdots,j_m\in [n]}(V_{S_k}^T)_{i_1j_1}a_{j_1j_2\cdots j_m}(V_{S_k})_{j_2i_2}
\cdots (V_{S_k})_{j_mi_m} \\
=&\sum_{j_1,j_2,\cdots,j_m\in [n]}(V_{S_k})_{j_1i_1}a_{j_1j_2\cdots j_m}(V_{S_k})_{j_2i_2}
\cdots (V_{S_k})_{j_mi_m}\\
=&\sum_{j_1,j_2,\cdots,j_m\in [n]}a_{j_1j_2\cdots j_m} ({\bf u}^k_{i_1})_{j_1}({\bf u}^k_{i_2})_{j_2}\cdots ({\bf u}^k_{i_m})_{j_m}\\
=&\langle\mathcal{A}, {\bf u}^k_{i_1} \circ{\bf u}^k_{i_2}\circ \cdots \circ {\bf u}^k_{i_m} \rangle.
\end{aligned}
\end{equation}
For any ${\bf x}\in S_0$, it follows that ${\bf x}\in S_k$ for some $k\in [r]$ such that
$$
{\bf x}=\lambda_1 {\bf u}^k_1+\lambda_2 {\bf u}^k_2+\cdots+\lambda_n {\bf u}^k_n,~~\sum_{i=1}^n\lambda_i=1,~~\lambda_i\geq 0,~\forall\, i\in [n].
$$
Thus,
$$
\begin{aligned}
\mathcal{A}{\bf x}^m=&\langle\mathcal{A}, (\lambda_1 {\bf u}^k_1+\lambda_2 {\bf u}^k_2+\cdots+\lambda_n {\bf u}^k_n)^m\rangle \\
=&\sum_{i_1,i_2,\cdots,i_m\in [n]}\lambda_{i_1}\lambda_{i_2}\cdots \lambda_{i_m}\langle\mathcal{A}, {\bf u}^k_{i_1} \circ{\bf u}^k_{i_2}\circ \cdots \circ {\bf u}^k_{i_m} \rangle \\
=&\sum_{i_1,i_2,\cdots,i_m\in [n]}(V_{S_k}^T\mathcal{A}V_{S_k})_{i_1i_2\cdots i_m}\lambda_{i_1}\lambda_{i_2}\cdots \lambda_{i_m} \\
=&(V_{S_k}^T\mathcal{A}V_{S_k}){\bf \lambda}^m,
\end{aligned}
$$
where ${\bf \lambda}=(\lambda_1,\lambda_2,\cdots,\lambda_n)\in \mathbb{R}_+^n$, and the third equality is obtained from (\ref{e43}).
By conditions and ${\bf \lambda}\in \mathbb{R}^n_+$, it holds that $\mathcal{A}{\bf x}^m\geq 0$ holds for any ${\bf x}\in S_0$; and hence, the desired results follow. \qed

To show the simplicial partition is fine enough, we will give a necessary condition for strictly copositivity of tensor.
For the standard simplex $S_0$ with a simplicial partition $\widetilde{S}=\{S_1,S_2,\cdots,S_r\}$, the vertices of simplex $S_k$ are denoted by ${\bf u}^k_1,{\bf u}^k_2,\cdots,{\bf u}^k_n$ for any $k\in [r]$. Let $d(\widetilde{S})$ denote the maximum diameter of a simplex in $\widetilde{S}$:
$$
d(\widetilde{S})=\max_{k\in [r]}\max_{i,j\in [n]}\|{\bf u}^k_i-{\bf u}^k_j\|_2.
$$

\bt\label{theorem43} Let $\mathcal{A}\in S_{m,n}$ be a strictly copositive tensor. Then, there exists $\varepsilon>0$ such that for all
finite simplicial partitions $\widetilde{S}=\{S_1,S_2,\cdots,S_r\}$ of $S_0$ with $d(\widetilde{S})<\varepsilon$, it follows that
$$
\langle \mathcal{A},~{\bf u}^k_{i_1}\circ {\bf u}^k_{i_2}\circ \cdots\circ {\bf u}^k_{i_m} \rangle > 0
$$
for all $k\in [r], i_j\in [n],j\in [m]$, where ${\bf u}^k_1, {\bf u}^k_2,\cdots, {\bf u}^k_n$
are vertices of the simplex $S_k$.
\et
\proof First of all, we define the following function:
$$
f({\bf w})=\langle\mathcal{A},~\underbrace{{\bf x}\circ{\bf y}\circ\cdots \circ{\bf z}}_m \rangle,~~\forall~{\bf w}:=({\bf x},{\bf y},\cdots,{\bf z})\in \underbrace{S_0\times S_0\times \cdots \times S_0}_m\in \mathbb{R}^{mn}.
$$
The strictly copositivity of $\mathcal{A}$ implies that $f({\bf w}_{{\bf x}})>0$ for any ${\bf w}_{{\bf x}}:=({\bf x}, {\bf x}, \cdots, {\bf x})\in S_0\times S_0\times \cdots \times S_0$. By continuity, it follows that, for any ${\bf x}\in S_0$, there exists $\varepsilon_{\bf x}>0$ such that
$$
f({\bf w})>0~~\mbox{\rm for~all}~{\bf w}~\mbox{\rm satisfying}~\|{\bf w}-{\bf w}_{{\bf x}}\|_2\leq \varepsilon_{\bf x}.
$$
Let $\varepsilon=\min_{{\bf x}\in S_0} \varepsilon_{\bf x}>0$. Then, it follows from uniformly continuity of $f$ that
\begin{equation}\label{e44}
f({\bf w})>0~~\mbox{\rm for~all}~{\bf w}~\mbox{\rm satisfying}~\|{\bf w}-{\bf w}_{{\bf x}}\|_2\leq \varepsilon~\mbox{\rm for~all}~{\bf x}\in S_0.
\end{equation}
For any simplicial partition $\widetilde{S}=\{S_1,S_2,\cdots,S_r\}$ of $S_0$ with $d(\widetilde{S})<\frac{1}{m}\varepsilon$, it holds that
$$
\|{\bf u}^k_i-{\bf u}^k_j\|_2\leq \frac{1}{m}\varepsilon,\quad \forall\, k\in [r],~\forall i,j\in [n].
$$
Moreover, for any ${\bf x}\in S_k\subseteq \widetilde{S}$, it follows that, for any $i\in [n]$,
$$
\|{\bf u}^k_i-{\bf x}\|_2\leq \|{\bf u}^k_i-{\bf u}^k_j\|_2\leq d(\widetilde{S})<\frac{1}{m}\varepsilon\quad \mbox{\rm for~some~} j\in [n],
$$
which implies that
$$
\|({\bf u}^k_{i_1}, {\bf u}^k_{i_2},\cdots, {\bf u}^k_{i_m})-{\bf w}_{{\bf x}}\|_2\leq \varepsilon,~~\forall\, i_j\in [n],~\forall\, j\in [m].
$$
Combining this with (\ref{e44}), we obtain that
$$
f({\bf u}^k_{i_1}, {\bf u}^k_{i_2},\cdots, {\bf u}^k_{i_m})=\langle \mathcal{A},~{\bf u}^k_{i_1}\circ {\bf u}^k_{i_2}\circ \cdots\circ {\bf u}^k_{i_m} \rangle > 0
$$
holds for all $k\in [r], i_j\in [n],j\in [m]$; and hence, the desired results follow. \qed

The following lemma gives a detection criterion for the case of a tensor being not copositive.
\bt\label{theorem44}
Let $\mathcal{A}\in S_{m,n}$ be given. Then, the following two assertions are equivalent.

{\bf (i)} Tensor $\mathcal{A}$ is not copositive.

{\bf (ii)} There exists $\varepsilon>0$ such that, for all simplicial partition $\widetilde{S}=\{S_1,S_2,\cdots,S_r\}$ of $S_0$ with $d(\widetilde{S})<\varepsilon$, there are at least one $k\in [r]$ and one $i\in [n]$ satisfying $\mathcal{A}({\bf u}^k_i)^m<0.$
\et
\proof It is obvious that {\bf (ii)} $\Rightarrow$ {\bf (i)}. To prove the converse, we assume that $\mathcal{A}$ is not copositive.
Then, there is ${\bf x}\in S_0$ such that $\mathcal{A}{\bf x}^m<0$. By continuity, there exists $\varepsilon>0$ such that
\begin{equation}\label{e44-1}
\mathcal{A}{\bf y}^m<0\quad \mbox{\rm for~all}~y~\mbox{\rm satisfying}~\|{\bf y}-{\bf x}\|_2<\varepsilon.
\end{equation}
For any simplicial partition $\widetilde{S}=\{S_1,S_2,\cdots,S_r\}$ of $S_0$ with $d(\widetilde{S})<\varepsilon$, there is at least one $k\in [r]$
such that ${\bf x}\in S_k$ with $\|{\bf u}^k_i-{\bf x}\|_2\leq d(\widetilde{S})<\varepsilon$. Thus, it follows from (\ref{e44-1}) that  $\mathcal{A}({\bf u}^k_i)^m<0$, which implies that the desired results follow. \qed

\setcounter{equation}{0}
\section{Detection algorithm based on simplicial partition}
Based on the results obtained in the last section, we can develop an algorithm to verify whether a tensor is
copositive or not, which is stated as follows.
%. The main idea can be roughly described as follows: Start from $\widetilde{S}=\{S_0\}$, where $S_0$ is the standard simplex of $\mathbb{R}^n_+$; and then, check whether there is a vertex ${\bf u}$ satisfying $\mathcal{A}{\bf u}^m<0$, or whether the copositivity criterion of Theorem \ref{theorem41} is satisfied. If the two cases do not appear, divide $S_0$ into two simplices  refine the partition and repeat the above process. Formally, this procedure is stated in Algorithm 5.1.
\vspace{4mm}

\begin{tabular}{@{}l@{}}
\hline
 \multicolumn{1}{c}{\bf Algorithm 5.1. Test whether a given symmetric tensor is copositive or not} \\
\hline
{\bf Input:} $\mathcal{A}\in S_{m,n}$  \\
\qquad Set $\widetilde{S}:=\{S_0\}$, where $S_0=conv\{{\bf e_1},{\bf e_2},\cdots,{\bf e_n}\}$ is the standard simplex      \\
\qquad while $\widetilde{S}\neq \emptyset$ do         \\
\qquad\qquad choose $S=conv\{{\bf u}_1,{\bf u}_2,\cdots,{\bf u}_n\}\in \widetilde{S}$ \\
\qquad\qquad if there exists $i\in [n]$ such that $\mathcal{A}{\bf u}^m_i<0$, then     \\
\qquad\qquad\qquad return ``$\mathcal{A}$ is not copositive''       \\
\qquad\qquad else if $\langle \mathcal{A}, {\bf u}_{i_1}\circ {\bf u}_{i_2}\circ \cdots \circ {\bf u}_{i_m}\rangle \geq 0$ for all $i_1,i_2,\cdots,i_m\in [n]$, then  \\
\qquad\qquad\qquad   $\widetilde{S}=\widetilde{S} \backslash \{S\}$         \\
\qquad\qquad else \\
\qquad\qquad\qquad simplicial partition $S=S_1\bigcup S_2$; and set $\widetilde{S}:=\widetilde{S}\backslash \{S\}\bigcup \{S_1, S_2\}$ \\
\qquad\qquad end if \\
\qquad end while \\
\qquad return `` $\mathcal{A}$ is copositive.'' \\
{\bf Output:} ``$\mathcal{A}$ is copositive'' or ``$\mathcal{A}$ is not copositive''.\\
\hline
\end{tabular}
\vspace{4mm}

From Theorems \ref{theorem43} and \ref{theorem44}, it is easy to see that the following result holds.
\bt\label{theorem51} In Algorithm 5.1, if the input symmetric tensor $\mathcal{A}$ is strictly copositive or $\mathcal{A}$
is not copositive, then the method will terminate in finitely many iterations.
\et

Thus, it is clear that Algorithm 5.1 can capture all strictly copositive tensors and non-copositive tensors. Unfortunately, when $\mathcal{A}$ is copositive but not strictly copositive, it is possible that the partition procedure of the algorithm leads to $d(\widetilde{S})\rightarrow 0$; and in this case, the algorithm dose not stop in general. This case also exists for the matrix detecting process \cite{Bundfuss08,Sponsel12}.
The reason for this is the following result.
\begin{proposition}\label{prop51} Suppose $\mathcal{A}\in S_{m,n}$ is copositive.
Let $S=conv\{{\bf u}_1,{\bf u}_2,\cdots,{\bf u}_n\}$ be a simplex with $\mathcal{A}{\bf u}_i^m>0$ for all $i\in [n]$.
If there exists ${\bf x}\in S\backslash \{{\bf u}_1,{\bf u}_2,\cdots,{\bf u}_n\}$ such that $\mathcal{A}{\bf x}^m=0$,
then there are $i_1,i_2,\cdots,i_m\in [n]$ such that
$\langle\mathcal{A}, {\bf u}_{i_1}\circ {\bf u}_{i_2}\circ \cdots \circ{\bf u}_{i_m}\rangle <0$.
\end{proposition}
\proof We prove this proposition by contradiction. Assume that
$$
\langle\mathcal{A}, {\bf u}_{i_1}\circ {\bf u}_{i_2}\circ \cdots \circ{\bf u}_{i_m}\rangle \geq0
$$
for all $i_1,i_2,\cdots,i_m\in [n]$.
By conditions, there is ${\bf x}\in S\backslash \{{\bf u}_1,{\bf u}_2,\cdots,{\bf u}_n\}$ such that
$\mathcal{A}{\bf x}^m=0$. It follows that there exist $\lambda_1,\ldots,\lambda_n\in \mathbb{R}_+$ such that ${\bf x}=\sum_{i=1}^n\lambda_i{\bf u}_i$ and $\sum_{i=1}^n\lambda_i=1$. Thus,
$$
\begin{aligned} \mathcal{A}{\bf x}^m=&\langle \mathcal{A}, (\lambda_1{\bf u}_1+\cdots+\lambda_n{\bf u}_n)^m \rangle \\
=&\sum_{i_1,i_2,\cdots,i_m\in [n]}\lambda_{i_1}\lambda_{i_2}\cdots \lambda_{i_m}\langle\mathcal{A}, {\bf u}_{i_1} \circ{\bf u}_{i_2}\circ \cdots \circ {\bf u}_{i_m} \rangle \\
\geq &\sum_{i=1}^n\lambda_i\mathcal{A}{\bf u}_i^m >0,
\end{aligned}
$$
which contradicts the result $\mathcal{A}{\bf x}^m=0$. Therefore, the desired results follow. \qed

To get rid of the termination problem with the given tensor $\mathcal{A}$ being copositive but not strictly copositive,
we can first try to check the copositivity of $\mathcal{A}$ by Algorithm 5.1. If it terminates in finitely many iterations, then
we get a correct answer; and if not, we can consider a relaxation form of $\mathcal{A}$, i.e.,
$\mathcal{B}=\mathcal{A}+\sigma\mathcal{E}$, where $\mathcal{E}$ is the tensor of all ones and $\sigma>0$ is a small tolerance.
\bt\label{theorem52} Let $\mathcal{A}\in S_{m,n}$ be given. Then, $\mathcal{A}$ is copositive if and only if
$\mathcal{B}=\mathcal{A}+\sigma\mathcal{E}$ is strictly copositive for any $\sigma>0$.
\et
\proof It is obvious that the necessary condition holds. For the sufficient statement, since
\begin{equation*}%\label{e51}
\mathcal{B}{\bf x}^m=(\mathcal{A}+\sigma\mathcal{E}){\bf x}^m=\mathcal{A}{\bf x}^m+\sigma>0,
\end{equation*}
for all ${\bf x}\in S_0$ and $\sigma>0$, by letting $\sigma\rightarrow 0$, we can obtain the desired result. \qed

From Theorems \ref{theorem51} and \ref{theorem52}, the following conclusion holds.
\bc\label{corol51} If the given tensor $\mathcal{A}$ is copositive but not strictly copositive, by replacing $\mathcal{A}$ by $\mathcal{B}=\mathcal{A}+\sigma\mathcal{E}$ for some $\sigma>0$, then Algorithm 5.1 terminates in finitely many iterations.
\ec

For $\sigma>0$, we call the symmetric tensor $\mathcal{A}$ a $\sigma$-copositive tensor with respect to simplex $S_0$, if $\mathcal{A}{\bf x}^m\geq -\sigma$ for all ${\bf x}\in S_0$.
And we have the following conclusion.
\begin{proposition}\label{proposition52}
Suppose $\mathcal{A}\in S_{m,n}$ is given. Let $S=conv\{{\bf u}_1,{\bf u}_2,\cdots,{\bf u}_n\}$ be a simplex and $\sigma>0$. If $\langle\mathcal{A}, {\bf u}_{i_1}\circ {\bf u}_{i_2}\circ \cdots \circ{\bf u}_{i_m}\rangle\geq -\sigma$ for all $i_1,i_2,\cdots,i_m\in [n]$, then $\mathcal{A}{\bf x}^m\geq -\sigma$ for any ${\bf x}\in S$.
\end{proposition}
\proof For any ${\bf x}\in S$, there exist $\lambda_1,\ldots,\lambda_n\in \mathbb{R}_+$ such that ${\bf x}=\sum_{i=1}^n\lambda_i{\bf u}_i$ and $\sum_{i=1}^n\lambda_i=1$. Thus,
$$
\begin{aligned} \mathcal{A}{\bf x}^m=&\langle \mathcal{A}, (\lambda_1{\bf u}_1+\cdots+\lambda_n{\bf u}_n)^m \rangle \\
=&\sum_{i_1,i_2,\cdots,i_m\in [n]}\lambda_{i_1}\lambda_{i_2}\cdots \lambda_{i_m}\langle\mathcal{A}, {\bf u}_{i_1} \circ{\bf u}_{i_2}\circ \cdots \circ {\bf u}_{i_m} \rangle \\
\geq& -\sigma,
\end{aligned}
$$
which implies that the desired result holds.\qed

\setcounter{equation}{0}
\section{Numerical examples}

In order to implement Algorithm 5.1, we specify Algorithm 5.1 as follows.
\vspace{4mm}

\begin{tabular}{@{}l@{}}
\hline
 \multicolumn{1}{c}{\bf Algorithm 6.1. Test whether a given symmetric tensor is copositive or not} \\
\hline
{\bf Input:} $\mathcal{A}\in S_{m,n}$  \\
\qquad Set $\widetilde{S}:=\{S_1\}$, where $S_1=conv\{{\bf e_1},{\bf e_2},\cdots,{\bf e_n}\}$ is the standard simplex      \\
\qquad Set $k:=1$ \\
\qquad while $k\neq 0$ do         \\
\qquad\qquad set $S:=S_k=conv\{{\bf u}_1,{\bf u}_2,\cdots,{\bf u}_n\}\in \widetilde{S}$ \\
\qquad\qquad if there exists $i\in [n]$ such that $\mathcal{A}{\bf u}^m_i<0$, then     \\
\qquad\qquad\qquad return ``$\mathcal{A}$ is not copositive''       \\
\qquad\qquad else if $\langle \mathcal{A}, {\bf u}_{i_1}\circ {\bf u}_{i_2}\circ \cdots \circ {\bf u}_{i_m}\rangle \geq 0$ for all $i_1,i_2,\cdots,i_m\in [n]$, then  \\
\qquad\qquad\qquad   set $\widetilde{S}:=\widetilde{S} \backslash \{S_k\}$ and $k:=k-1$    \\
\qquad\qquad else \\
\qquad\qquad\qquad set  \\
\qquad\qquad\qquad\qquad $S_k:=conv\{{\bf u}_1,\cdots,{\bf u}_{p-1},{\bf v},{\bf u}_{p+1},\cdots,{\bf u}_n\}$;  \\
\qquad\qquad\qquad\qquad $S_{k+1}:=conv\{{\bf u}_1,\cdots,{\bf u}_{q-1},{\bf v},{\bf u}_{q+1},\cdots,{\bf u}_n\}$, \\
\qquad\qquad\qquad\qquad where ${\bf v}=\frac{{\bf u_p}+{\bf u_q}}{2},\, [p,q]=arg\max_{i,j\in [n]}\|{\bf u_i}-{\bf u_j}\|_2$ and $p<q$.\\
\qquad\qquad\qquad set $\widetilde{S}:=\widetilde{S}\backslash \{S\}\bigcup \{S_k, S_{k+1}\}$ and $k:=k+1$ \\
\qquad\qquad end if \\
\qquad end while \\
\qquad return `` $\mathcal{A}$ is copositive.'' \\
{\bf Output:} ``$\mathcal{A}$ is copositive'' or ``$\mathcal{A}$ is not copositive''.\\
\hline
\end{tabular}
\vspace{4mm}

In this section, we use this specified version of Algorithm 5.1 to detect whether a tensor is copositive or not. All experiments are finished in Matlab2014b on a Philips desktop computer with Intel(R) Core(TM)2 Duo CPU E8500 @ 3.16GHz 3.17 GHz and 4 GB of RAM. We detect several classes of tensors from three aspects, which are given in the following three parts, respectively.

{\bf Part 1}. Suppose that $\mathcal{B}\in S_{m,n}$ is a nonnegative tensor and $\rho(\mathcal{B})$ denotes its spectral radius. Let $\mathcal{I}\in S_{m,n}$ denote the identity tensor. Then, by Definition 3.1 and Theorem 3.12 of \cite{Zhang12}, we have the following results: The tensor $\eta \mathcal{I}-\mathcal{B}$ is copositive if and only if $\eta\geq \rho(\mathcal{B})$; and the tensor $\eta \mathcal{I}-\mathcal{B}$ is strictly copositive if and only if $\eta> \rho(\mathcal{B})$. Based on these results, we construct several tensors for testing. We first test the following specific tensors.
\begin{example}\label{example1}
Suppose that $\mathcal{A}\in S_{3,3}$ (or $\mathcal{A}\in S_{4,4}$) is given by
\begin{eqnarray}\label{E-exam1-1}
\mathcal{A}=\eta \mathcal{I}-\mathcal{B},
\end{eqnarray}
where $\mathcal{B}\in S_{3,3}$ (or $\mathcal{B}\in S_{4,4}$) is a tensor of ones and $\eta$ is specified in our numerical results.
\end{example}

The numerical results are given in Table 1, where ``$\rho$" denotes the spectral radius of the tested tensor, ``IT" denotes the number of iterations, ``CPU(s)" denotes the CPU time in seconds, and ``Result" denotes the output result in which ``No" denotes the output result that the tested tensor is not copositive and ``Yes" denotes the output result that the tested tensor is copositive.

\begin{table}[ht] %\label{table1}
  \caption{The numerical results of the problem in Example \ref{example1}}
  \begin{center}
    \begin{tabular}[c]
      {| c | c | c | c | c | c | c|}
      \hline
      $m$  & $n$   &$\rho$ & $\eta$ & IT      & CPU(s) & Result\\
      \hline
      \hline
           &       &       & 1      & 2       &0.078     & No   \\  \cline{4-7}
           &       &       & 8.99   & 43      &0.437     & No   \\ \cline{4-7}
        3  & 3     & 9     & 9      & $>100$  &          &      \\  \cline{4-7}
           &       &       & 9.01   & 59      &0.593     & Yes   \\  \cline{4-7}
           &       &       & 19     & 11      &0.172     & Yes   \\ \hline
           &       &       & 10     & 14      &0.499     & No   \\ \cline{4-7}
        4  & 4     & 64    & 64     & 63      &2.32      & Yes   \\ \cline{4-7}
           &       &       & 74     & 63      &2.14      & Yes   \\ \hline
    \end{tabular}
  \end{center}
\end{table}

Next, we test some randomly generated tensors with the form being same as the one by (\ref{E-exam1-1}).

\begin{example}\label{example2}
Suppose that $\mathcal{A}\in S_{m,n}$ is given by (\ref{E-exam1-1}), where $\mathcal{B}\in S_{m,n}$ is randomly generated with all its elements are in the interval $(0,1)$.
\end{example}

In our experiments, we use the higher order power method to compute the spectral radius $\rho$ of every tensor $\mathcal{B}$. $m$, $n$ and $\eta$ are specified in our numerical results. For the same $m$ and $n$, we generate randomly every tested problem 10 times, the numerical results are shown in Table 2, where ``MinIT" and ``MaxIT" denote the minimal number and the maximal number of iterations among ten times for every tested problem, respectively, ``MinCPU(s)" and ``MaxCPU(s)" denote the smallest and the largest CPU times in second among ten times for every tested problem, respectively, ``Nyes" denotes the number of the output result that the tested tensor is copositive, and ``Nno" denotes the number of the output result that the tested tensor is not copositive. The same notations are also used in Table 6.3.

\begin{table}[ht] %\label{table1}
  \caption{The numerical results of the problem in Example \ref{example2}}
  \begin{center}
    \begin{tabular}[c]
      {| c | c | c | c | c | c | c| c | c |}
      \hline
      $m$  & $n$   & $\eta$       & MinIT   & MaxIT       &MinCPU(s) &MinCPU(s) & Nyes  & Nno\\
      \hline
      \hline
        3  & 3     & $\rho-1$  & 6       & 25          &0.0624    &0.234     &       & 10    \\  \cline{3-9}
           &       & $\rho+1$  & 19      & 19          &0.187     &0.187     & 10    &       \\  \cline{3-9}
           &       & $\rho+10$ & 11      & 11          &0.109     &0.125     & 10    &      \\  \hline
        3  & 4     & $\rho-1$  & 21      & 65          &0.39      &1.19      &       & 10    \\  \cline{3-9}
           &       & $\rho+1$  & 63      & 75          &1.14      &1.36      & 10    &       \\  \cline{3-9}
           &       & $\rho+10$ & 49      & 53          &0.905     &0.983     & 10    &      \\  \hline
        4  & 3     & $\rho-1$  & 17      & 17          &0.234     &0.312     &       & 10    \\  \cline{3-9}
           &       & $\rho+1$  & 27      & 31          &0.39      &0.452     & 10    &       \\  \cline{3-9}
           &       & $\rho+10$ & 19      & 19          &0.265     &0.296     & 10    &      \\  \hline
        4  & 4     & $\rho-1$  & 21      & 25          &0.686     &0.827     &       & 10    \\  \cline{3-9}
           &       & $\rho+1$  & 65      & 91          &2.11      &3.15      & 10    &       \\  \cline{3-9}
           &       & $\rho+10$ & 63      & 63          &2.09      &2.22      & 10    &      \\  \hline
        6  & 3     & $\rho-1$  & 20      & 28          &0.562     &0.811     &       & 10    \\  \cline{3-9}
           &       & $\rho+1$  & 43      & 47          &1.25      &1.36      & 10    &       \\  \cline{3-9}
           &       & $\rho+10$ & 27      & 27          &0.764     &0.796     & 10    &      \\  \hline
    \end{tabular}
  \end{center}
\end{table}

{\bf Part 2}. It is obvious that any nonnegative tensor is copositive. By Corollary 6.1 of \cite{Song15}, we also know that for any $\mathcal{A}\in S_{m,n}$, if $\mathcal{A}$ is (strictly) copositive, then ($a_{ii\cdots i}>0$) $a_{ii\cdots i}\geq 0$ for all $i\in [n]$. Based on these results, we consider to detect the tensors given in the following example.

\begin{example}\label{example3}
(i) Consider the tensor $\mathcal{A}\in S_{m,n}$ which is randomly generated with all its elements are in the interval $(0,1)$; (ii) we set $\mathcal{B}:=\mathcal{A}$ and $b_{11\cdots 1}=-1$.
\end{example}

In our experiments, for the same $m$ and $n$, we generate randomly every tested problem 10 times, the numerical results are shown in Table 3.

\begin{table}[ht] %\label{table1}
  \caption{The numerical results of the problem in Example \ref{example3}}
  \begin{center}
    \begin{tabular}[c]
      {| c | c | c | c | c | c | c| c | c |}
      \hline
      $m$  & $n$   & Tensor        &MinIT  & MaxIT       &MinCPU(s) &MinCPU(s) & Nyes  & Nno\\
      \hline
      \hline
        3  & 3     & $\mathcal{A}$ & 1     & 1          &0.0624    &0.0936     & 10    &      \\  \cline{3-9}
           &       & $\mathcal{B}$ & 1     & 1          &0.0624    &0.078      &       & 10   \\  \hline
        3  & 4     & $\mathcal{A}$ & 1     & 1          &0.0624    &0.078      & 10    &      \\  \cline{3-9}
           &       & $\mathcal{B}$ & 1     & 1          &0.0624    &0.078      &       & 10   \\  \hline
        4  & 3     & $\mathcal{A}$ & 1     & 1          &0.078     &0.078      & 10    &      \\  \cline{3-9}
           &       & $\mathcal{B}$ & 1     & 1          &0.0624    &0.078      &       & 10   \\  \hline
        4  & 4     & $\mathcal{A}$ & 1     & 1          &0.0936    &0.125      & 10    &      \\  \cline{3-9}
           &       & $\mathcal{B}$ & 1     & 1          &0.0624    &0.078      &       & 10   \\  \hline
        6  & 3     & $\mathcal{A}$ & 1     & 1          &0.0936    &0.125      & 10    &      \\  \cline{3-9}
           &       & $\mathcal{B}$ & 1     & 1          &0.0936    &0.109      &       & 10   \\  \hline
    \end{tabular}
  \end{center}
\end{table}

{\bf Part 3}. It is well known that there is a one-to-one relationship between the homogeneous polynomial and the symmetric tensor. In this part, we consider several tensors which come from several famous homogeneous polynomials. For convenience, we use the following notation: for any integers $i_1,i_2,\ldots,i_m$, we use $\pi(i_1i_2\cdots i_m)$ to denote a permutation of $i_1i_2\cdots i_m$, and $S_{\pi(i_1i_2\cdots i_m)}$ to denote the set of all these permutations.
\begin{example}\label{example4}
Suppose that $\mathcal{A}\in S_{6,3}$ is given by
\begin{eqnarray*}
\left\{\begin{array}{ll}
\sum_{i_1i_2i_3i_4i_5i_6\in S_{\pi(111122)}}a_{i_1i_2i_3i_4i_5i_6}=1,\\
\sum_{i_1i_2i_3i_4i_5i_6\in S_{\pi(112222)}}a_{i_1i_2i_3i_4i_5i_6}=1,\\
a_{333333}=1,\\
\sum_{i_1i_2i_3i_4i_5i_6\in S_{\pi(112233)}}a_{i_1i_2i_3i_4i_5i_6}=-3,
\end{array}\right.
\end{eqnarray*}
\end{example}

The corresponding polynomial of the tensor $\mathcal{A}$ given in this example is
$$
f(x,y,z)=x^4y^2+x^2y^4+z^6-3x^2y^2z^2.
$$
This is the famous Motzkin polynomial, which is non-negative but not a sum of squares; and hence, the tensor $\mathcal{A}$ is copositive. It is easy to see that the tensor $\mathcal{A}$ is not strictly copositive. We use Algorithm 5.1 to test this tensor, the algorithm does not terminate within $100$ iterations. We use Algorithm 5.1 to test the tensor $\mathcal{A}+\sigma\mathcal{E}$ with $\sigma>0$, however, the algorithm can correctly detect the copositivity of the tensor; and we list several cases as follows.
\begin{itemize}
\item When $\sigma=0.01$, the algorithm can correctly detect the copositivity of the tensor with $11$ iterations in $0.406$ seconds;
\item when $\sigma=0.001$, the algorithm can correctly detect the copositivity of the tensor with $27$ iterations in $0.874$ seconds; and
\item when $\sigma=0.0001$, the algorithm can correctly detect the copositivity of the tensor with $71$ iterations in $2.25$ seconds.
\end{itemize}

\begin{example}\label{example5}
Suppose that $\mathcal{A}\in S_{6,3}$ is given by
\begin{eqnarray*}
\left\{\begin{array}{ll}
a_{111111}=1,\;\; a_{222222}=1,\;\; a_{333333}=1,\\
\sum_{i_1i_2i_3i_4i_5i_6\in S_{\pi(111122)}}a_{i_1i_2i_3i_4i_5i_6}=-1,\\
\sum_{i_1i_2i_3i_4i_5i_6\in S_{\pi(112222)}}a_{i_1i_2i_3i_4i_5i_6}=-1,\\
\sum_{i_1i_2i_3i_4i_5i_6\in S_{\pi(111133)}}a_{i_1i_2i_3i_4i_5i_6}=-1,\\
\sum_{i_1i_2i_3i_4i_5i_6\in S_{\pi(113333)}}a_{i_1i_2i_3i_4i_5i_6}=-1,\\
\sum_{i_1i_2i_3i_4i_5i_6\in S_{\pi(222233)}}a_{i_1i_2i_3i_4i_5i_6}=-1,\\
\sum_{i_1i_2i_3i_4i_5i_6\in S_{\pi(223333)}}a_{i_1i_2i_3i_4i_5i_6}=-1,\\
\sum_{i_1i_2i_3i_4i_5i_6\in S_{\pi(112233)}}a_{i_1i_2i_3i_4i_5i_6}=3,
\end{array}\right.
\end{eqnarray*}
\end{example}

The corresponding polynomial of the tensor $\mathcal{A}$ given in this example is
$$
f(x,y,z)=x^6+y^6+z^6-x^4y^2-x^2y^4-x^4z^2-x^2z^4-y^4z^2-y^2z^4+3x^2y^2z^2.
$$
This is the famous Robinson polynomial, which is non-negative but not a sum of squares; and hence, the tensor $\mathcal{A}$ is copositive. It is easy to see that the tensor $\mathcal{A}$ is not strictly copositive. We use Algorithm 5.1 to test this tensor, the algorithm does not terminate within $100$ iterations. We use Algorithm 5.1 to test the tensor $\mathcal{A}+\sigma\mathcal{E}$ with $\sigma>0$, however, the algorithm can correctly detect the copositivity of the tensor; and we list several cases as follows.
\begin{itemize}
\item When $\sigma=0.01$, the algorithm can correctly detect the copositivity of the tensor with $11$ iterations in $0.406$ seconds;
\item when $\sigma=0.001$, the algorithm can correctly detect the copositivity of the tensor with $27$ iterations in $0.842$ seconds; and
\item when $\sigma=0.0001$, the algorithm can correctly detect the copositivity of the tensor with $83$ iterations in $2.5$ seconds.
\end{itemize}

\begin{example}\label{example6}
Suppose that $\mathcal{A}\in S_{6,3}$ is given by
\begin{eqnarray*}
\left\{\begin{array}{ll}
\sum_{i_1i_2i_3i_4i_5i_6\in S_{\pi(111122)}}a_{i_1i_2i_3i_4i_5i_6}=1,\\
\sum_{i_1i_2i_3i_4i_5i_6\in S_{\pi(222233)}}a_{i_1i_2i_3i_4i_5i_6}=1,\\
\sum_{i_1i_2i_3i_4i_5i_6\in S_{\pi(333311)}}a_{i_1i_2i_3i_4i_5i_6}=1,\\
\sum_{i_1i_2i_3i_4i_5i_6\in S_{\pi(112233)}}a_{i_1i_2i_3i_4i_5i_6}=-3,
\end{array}\right.
\end{eqnarray*}
\end{example}

The corresponding polynomial of the tensor $\mathcal{A}$ given in this example is
$$
f(x,y,z)=x^4y^2+y^4z^2+z^4x^2-3x^2y^2z^2.
$$
This is the famous Choi-Lam polynomial, which is non-negative but not a sum of squares; and hence, the tensor $\mathcal{A}$ is copositive. It is easy to see that the tensor $\mathcal{A}$ is not strictly copositive. We use Algorithm 5.1 to test this tensor, the algorithm does not terminate within $100$ iterations. We use Algorithm 5.1 to test the tensor $\mathcal{A}+\sigma\mathcal{E}$ with $\sigma>0$, however, the algorithm can correctly detect the copositivity of the tensor; and we list several cases as follows.
\begin{itemize}
\item When $\sigma=0.01$, the algorithm can correctly detect the copositivity of the tensor with $5$ iterations in $0.218$ seconds;
\item when $\sigma=0.001$, the algorithm can correctly detect the copositivity of the tensor with $27$ iterations in $0.858$ seconds; and
\item when $\sigma=0.0001$, the algorithm can correctly detect the copositivity of the tensor with $41$ iterations in $1.29$ seconds.
\end{itemize}

From the numerical results given in Part 1-Part 3, we can see that Algorithm 6.1 is effective for the problems we tested.

\setcounter{equation}{0}
\section{Conclusions}

In this paper, we proposed new criteria to judge whether a tensor is (strictly) copositive or not, including three necessary conditions which are given based on binomial expansion, principle subtensor and convex combination, respectively; and several necessary conditions or sufficient conditions which are investigated by taking advantage of the simplicial partition. These theoretical results can be viewed as extensions of those obtained in the case of matrix. Moreover, by the obtained criteria based on the simplicial partition, we proposed a detection algorithm for the copositive tensor. The preliminary numerical results demonstrate that the proposed algorithm is effective.

\end{document}